\documentclass[letterpaper, 10 pt, conference]{ieeeconf}
\IEEEoverridecommandlockouts                              

\usepackage{amsfonts}
\usepackage{amsmath,amssymb}         
\usepackage{hyperref,cite}
\usepackage{cleveref}
\usepackage{color}
\usepackage{algorithm}
\usepackage{algorithmic}
\let\labelindent\relax
\usepackage{enumitem}
\usepackage{float}
\usepackage{graphicx}
\usepackage{booktabs}
\usepackage{threeparttable}
\usepackage{subfigure}
\usepackage{multirow}

\usepackage{etoolbox}
\makeatletter
\patchcmd{\@makecaption}
  {\scshape}
  {}
  {}
  {}
\makeatother

\newtheorem{assumption}{Assumption}
\DeclareMathOperator*{\argmin}{argmin} 
\newcommand{\diag}{\mathop{\rm diag}\nolimits}

\begin{document}
\title{\LARGE \bf
A construction-free coordinate-descent augmented-Lagrangian method
for embedded linear MPC based on ARX models
}

\author{Liang Wu$^{1}$, Alberto Bemporad$^{1}$%
\thanks{The authors are with the IMT School for Advanced Studies Lucca, Italy,
        {\tt\small \{liang.wu,alberto.bemporad\}@imtlucca.it}}%
}
\thispagestyle{empty}
\pagestyle{empty}
\maketitle

\begin{abstract}
This paper proposes a construction-free algorithm for solving linear MPC problems based on autoregressive with exogenous terms (ARX) input-output models. The solution algorithm relies on a coordinate-descent augmented Lagrangian (CDAL) method previously proposed by the authors, which we adapt here to exploit the special structure of ARX-based MPC. The CDAL-ARX algorithm enjoys the construction-free feature, in that it avoids explicitly constructing the quadratic programming (QP) problem associated with MPC, which would eliminate construction cost when the ARX model changes/adapts online. For example, the ARX model parameters are dependent on linear parameter-varying (LPV) scheduling signals, or recursively adapted from streaming input-output data with cheap computation cost, which make the ARX model widely used in adaptive control. Moreover, the implementation of the resulting CDAL-ARX algorithm is matrix-free and library-free, and hence amenable for deployment in industrial embedded platforms. We show the efficiency of CDAL-ARX in two numerical examples, also in comparison with MPC implementations based on other general-purpose quadratic programming solvers.
\end{abstract}


\section{Introduction}
Model Predictive Control (MPC) is an advanced technique to control multi-input multi-output systems subject to constraints, and its core idea to predict the evolution of the controlled system by means of a dynamical model, solve an optimization problem over a finite time horizon, only implement the control input at the current time, and then repeat the optimization again at the next sample~\cite{qin2003survey}.
In earlier MPC developments, some methods shared the same receding horizon control idea, under different names. The Model Predictive Heuristic Control (MPHC) \cite{richalet1978model}, the Model Algorithmic Control (MAC) \cite{ROUHANI1982401} used a finite impulse response model, the Dynamic Matrix Control (DMC) employed a truncated step-response model \cite{cutler1980dynamic}, and the Generalized Predictive Control (GPC) involved a transfer function model \cite{CLARKE1987137}. As the MPC field has grown, state-space (SS) models replaced input-output (I/O) models, and most MPC theory is based on SS formulations~\cite{mayne2014model}. 

However, in industrial control applications, MPC based on input-output models, such as the autoregressive model with exogenous terms (ARX) model, may still be preferable \cite{qin2003survey}, for two main reasons: (1) there is no need of a state-observer; (2) I/O models are easier to identify and to adapt online (such as using recursive least-squares or Kalman Filter algorithms), which makes them widely used in adaptive control \cite{aastrom2013adaptive}.
In particular, the latter is particularly appealing in practical cases in which the dynamics
of the systems changes during operations, such as in the case of changes of mass and inertia in rockets due to fuel consumption, wear of heating equipment in chemical processes, and many others. In fact, an observable SS model can be equivalently transformed into an ARX model and our previous work \cite{wu2022equivalence} show the equivalence of SS-based MPC and ARX-based MPC problems.
In \cite{wu2022equivalence}, we proposed an alternative for the acquisition of ARX model based on the first-principle-based modeling paradigm, rather than the data-driven identification paradigm. This allows us to acquire ARX models using many existing first-principles based models in different engineering fields. The resulted interpretative ARX model can be adopted in adaptive MPC framework by adding an online updating scheme for the ARX model, see our work \cite{wu2022interpretative}.

A common practice in MPC is to first formulate an quadratic programming (QP) problem in terms of a control-oriented prediction model and MPC parameters, and then pass it to the optimization solver. Such a problem construction step can be performed offline when the prediction model is fixed, otherwise it requires to be repeated online when the prediction model or MPC parameters are varying. In such varying cases, the online computation time includes both constructing and solving the QP problem associated with MPC. Indeed, often constructing and solving the MPC problem have comparable costs, such as when warm-starting strategies are employed and set-points change slowly. The online construction of the MPC problem becomes necessary in the adaptive ARX-based MPC framework \cite{wu2022interpretative}, and the case of linear parameter varying ARX (LPV-ARX) models, in which model parameters depend on a measured time-varying signal, the so-called scheduling variable. Thus, an construction-free ARX-MPC algorithm, in that MPC-to-QP construction is explicitly eliminated, would significantly save the computational loads in those cases.

\subsection{Related works and Contribution}
Some ARX-based MPC algorithms in the literature first convert the ARX model into SS form, treat the problem as a standard SS-based MPC problem \cite{huusom2010arx}, and then construct and solve a condensed or sparse quadratic programming (QP) problem. In fact, the ARX-to-SS transformation is not necessary in condensed and sparse MPC-to-QP constructions, which only depend on whether eliminating or keeping the ARX output variables. Choosing the condensed or sparse construction only depends on the total online computation cost (constructing and solving) \cite{kouzoupis2015towards}, in time-varying ARX-based MPC problems. The OSQP solver, based on the alternating direction method of multipliers (ADMM), can directly consume the ARX model as equality constraints, which is the sparse QP formulation by keeping the output variables of the ARX model. However, it still employs an explicit MPC-to-QP construction to formulate the equality constraint matrix, and more importantly, the OSQP solver needs to repeatedly factorize and cache the Hessian matrix of the quadratic objective at each sampling time, in time-varying ARX-based MPC problems. In \cite{saraf2017fast}, the dynamic equality constraint from the ARX model was relaxed by using a large penalty parameter, and it resulted in an ill-conditioning bounded variable least-squares (BVLS) problem, although the active-set based method was used to mitigate the numerical difficulties to some extent. Besides computation efficiency, easy-to-deployment of an ARX-MPC algorithm should also be considered, in which the code simplicity and library-dependency are important. In this respect, compared to the active-set or interior-point based methods, the first-order method, such as the primal or dual fast gradient method, the ADMM method, is simpler but also becomes complicated in time-varying cases, in that some offline operations have to be performed online \cite{kouzoupis2015first}.

This paper proposes a simple and efficient algorithm for solving ARX-based MPC problems. Based on the coordinate-descent augmented Lagrangian method, the resulting CDAL-ARX algorithm enjoys three mian features: ($i$) it is \emph{construction-free}, in that it avoids the online MPC-to-QP construction in time-varying ARX cases to save computation cost; ($ii$) is \emph{matrix-free}, in that it avoids multiplications and factorizations of matrices, which are required by other first-order methods in time-varying ARX cases; and ($iii$) is \emph{library-free}, as our 150-lines of C-code implementation is without any library dependency, which matters in embedded deployment.



\section{ARX-based MPC problem formulation}
\label{sec:ARX-based MPC}
Consider the multi-input multi-output (MIMO) ARX model described by
\begin{equation}\label{ARX}
y_{t}=\sum_{i=1}^{n_{a}} A(i) y_{t-i}+\sum_{i=1}^{n_{b}} B(i) u_{t-i}
\end{equation}
where $y_{t}\in\mathbb{R}^{n_y}$ and $u_{t}\in \mathbb{R}^{n_u}$ are the output and input of the system, respectively, $A(i)\in \mathbb{R}^{n_y \times n_y}$,
$i=1,\ldots, n_a$, and $B(i)\in \mathbb{R}^{n_y \times n_u}$,
$i=1,\ldots, n_b$, and $n_a$, $n_b$ define the model order of the ARX model.

This paper considers the following MPC tracking formulation based the ARX model~\eqref{ARX}
\begin{equation}\label{ARX_MPC}
\begin{array}{ll}
\min_{Y,U,\Delta U} & \displaystyle{\frac{1}{2}\sum_{t=1}^{T}\left\|\left(y_{t}-r_{t}\right)\right\|_{W^{y}}^{2}+\left\| \Delta u_{t-1}\right\|_{W^{\Delta u}}^{2}}\\
\mbox{s.t.} & \hspace*{-2em}y_{t}=\displaystyle{\sum_{i=1}^{n_{a}} A(i) y_{t-i}+\sum_{i=1}^{n_{b}} B(i) u_{t-i},\ t=1,\ldots,T}\\
&\hspace*{-2em} \Delta u_{t}=u_{t}-u_{t-1},\ t=0,\ldots,T-1\\
&\hspace*{-2em} y_{\rm min } \leq y_{t} \leq y_{\rm max },\ t=1,\ldots,T\\
&\hspace*{-2em} u_{\rm min } \leq u_{t} \leq u_{\rm max },\ t=0,\ldots,T-1\\
&\hspace*{-2em} \Delta u_{\rm min } \leq \Delta u_{t} \leq \Delta u_{\rm max }, t=0,\ldots,T-1
\end{array}
\end{equation}
where $T$ is the prediction horizon, $W^y\succeq 0$ and $W^{\Delta u} \succeq 0$ are positive semi-definiteness diagonal matrices on the outputs and the input increments, respectively, $r_{t}$, $t=1,\ldots, T$ are the future desired set-point vectors, $\Delta u_{t-1}$ are the input increments, $[y_{\rm min},y_{\rm max}]$, $[u_{\rm min},u_{\rm max}]$, and $[\Delta u_{\rm min},\Delta u_{\rm max}]$ define box constraints on outputs, inputs, and input increments, respectively, and
$Y=(y_1,\ldots,y_T)$, $U=(u_0,\ldots,u_{T-1})$, and $\Delta U=(\Delta u_0,\ldots,\Delta u_{T-1})$
are the optimization variables.

\section{Coordinate Descent Augmented Lagrangian method}\label{sec:CDAL method}
In~\cite{wu2021simple}, we proposed a coordinate-descent augmented-Lagrangian (CDAL) method for SS-based MPC problems. We want to adapt here the method to solve problem~\eqref{ARX_MPC} without computing a state-space realization of the ARX model~\eqref{ARX}, while retaining the construction-free, matrix-free, and library-free properties of CDAL.

\subsection{Augmented Lagrangian method}
The following assumptions are needed to ensure the convergence of the Augmented Lagrangian method.
\begin{assumption}\label{assump1}
Problem~\eqref{ARX_MPC} has a feasible solution.
\end{assumption}
Note that Assumption~\eqref{assump1} is satisfied in all practical situations in which the reference $r_t$ is far enough from the output bounds and the prediction horizon $T$ is long enough.
\begin{assumption}\label{assump2}
The equality constraint matrix arising from stacking all the equality
constraints in~\eqref{ARX_MPC} is full rank at the optimal solution of 
the problem. 
\end{assumption}

Let $\mathcal{Y}$, $\mathcal{U}$, and $\Delta \mathcal{U}$ denote the hyper-boxes on $Y$, $U$, and $\Delta U$, respectively, defined by the box constraints in~\eqref{ARX_MPC}, respectively. The bound-constrained Augmented Lagrangian function $\mathcal{L}_{\rho}:  \mathcal{Y} \times \mathcal{U} \times \Delta\mathcal{U} \times \mathbb{R}^{Tn_y} \times \mathbb{R}^{T n_u} \rightarrow \mathbb{R}$ is given by
\begin{equation}
\begin{aligned}
\mathcal{L}_{\rho}&(Y,U,\Delta U,\Lambda,\Gamma)=\frac{1}{2} \sum_{t=1}^{T}\left\|\left(y_{t}-r_{t}\right)\right\|_{W y}^{2}+\left\|\Delta u_{t-1}\right\|_{W \Delta u}^{2} \\&
+\sum_{t=1}^{T} \lambda_{t}^{\prime}\left(\sum_{i=1}^{n_{a}} A(i) y_{t-i}+\sum_{i=1}^{n_{b}} B(i) u_{t-i}-y_{t}\right) \\
&+\sum_{t=1}^{T} \gamma_{t}^{\prime}\left(u_{t-2}+\Delta u_{t-1}-u_{t-1}\right) \\
&+\frac{\rho}{2} \sum_{t=1}^{T}\left\|\sum_{i=1}^{n_{a}} A(i) y_{t-i}+\sum_{i=1}^{n_{b}} B(i) u_{t-i}-y_{t}\right\|_{2}^{2} \\
&+\frac{\rho}{2} \sum_{t=1}^{T}\left\|u_{t-2}+\Delta u_{t-1}-u_{t-1}\right\| \|_{2}^{2}
\end{aligned}
\end{equation}
where $\Lambda = \{\lambda_{t}\}$ and $\Gamma = \{\gamma_{t}\}, \forall t=1,\ldots,T$ are the dual vectors associated to the equality constraints induced by the ARX model and the input increments, respectively, and $\rho$ is the penalty parameter. According to~\cite{bertsekas2014constrained}, the scaled AL method (ALM) iterates the following updates
\begin{subequations}\label{ALM}
\begin{eqnarray}
(Y^{k},U^{k},\Delta U^{k})= \argmin\frac{1}{\rho}\mathcal{L}_{\rho}(Y,U,\Delta U, \Lambda^{k-1},\Gamma^{k-1})\label{subALM}\\
\begin{aligned}
\lambda_{t}^{k}= \lambda_{t}^{k-1}+
\sum_{i=1}^{n_{a}} A(i) y_{t-i}^{k}+\sum_{j=1}^{n_{b}} B(i) u_{t-i}^{k}-y_{t}^{k}\\
,\forall t=1, \ldots, T
\end{aligned}\label{subALM_lambda}\\
\gamma_{t}^{k}=\gamma_{t}^{k-1}+u_{t-1}^{k}+\Delta u_{t}^{k}-u_{t}^{k}, \forall t=1,\ldots,T \label{subALM_beta}
\end{eqnarray}
\end{subequations}
The minimization step~\eqref{subALM} updates the primal vector, Steps~\eqref{subALM_lambda} and  \eqref{subALM_beta} update the dual vectors. We refer the reader to~\cite{bertsekas2014constrained}
for a well-known convergence proof of ALM under Assumptions~\ref{assump1},~\ref{assump2}. To improve the speed of convergence of ALM,~\cite{kang2015inexact} proposed an accelerated version of ALM whose convergence rate is $O(1/k^2)$ for linearly constrained convex programs by using Nesterov's acceleration technique~\cite{Nes83}. The accelerated ALM algorithm for MPC problems has been summarized in our previous work~\cite{wu2021simple}.

\subsection{Coordinate-descent method}
Sub-problem~\eqref{subALM} is a strongly convex box-constrained QP problem, which can be solved by many methods. Among others, as we showed in~\cite{wu2021simple}, problem~\eqref{subALM}
can be solved by a simple coordinate-descent method, which minimizes the objective function along only one coordinate direction at each iteration while keeping the other coordinates fixed~\cite{luo1992convergence}. A convergence proof of at-least linear convergence when solving convex differentiable minimization problems was shown in \cite{luo1992convergence}. Under Assumption~\ref{assump1} (non-emptyness of the feasible set) and since the objective function $\mathcal{L}_{\rho}(\cdot)$ is continuously differentiable and convex with respect to each coordinate, the CD method proceeds repeatedly for $k=1,2,\ldots,$ as follows:
\begin{subequations}\label{CD}
\begin{gather}
\text{choose}~j_k \in \left\{1,2,\ldots,n_z\right\} \\
z_{j_k}^{k} = \argmin_{z_{j_k} \in \mathcal{Z}} \frac{1}{\rho}\mathcal{L}_\rho(z_{j_k},z_{\neq j_k}^{k-1},\Lambda^{k-1},\Gamma^{k-1})
\end{gather}
\end{subequations}
where $z=\left[\begin{array}{ccccccc}
y_{1}^{\prime} & u_{0}^{\prime} & \Delta u_{0}^{\prime} & \ldots & y_{T}^{\prime} & u_{T-1}^{\prime} & \Delta u_{T-1}^{\prime}
\end{array}\right]^{\prime}$ is the optimization vector, 
$z\in\mathcal{Z}\triangleq\mathcal{Y} \times \mathcal{U} \times \Delta \mathcal{U}$,
$\mathcal{Z}\subseteq\mathbb{R}^{n_z}$, $n_z\triangleq T(n_y+n_u+n_u)$. We denote  by $\mathcal{L}_\rho(z_{j_k},z_{\neq j_k}^{k-1},\Lambda^{k-1},\Gamma^{k-1})$
the value $\mathcal{L}_\rho(z,\Lambda^{k-1},\Gamma^{k-1})$ when $z_{\neq j_k}=z_{\neq j_k}^k$ is fixed.
Here $z_{\neq j_k}$ denotes the subvector obtained from $z$
by eliminating its $j_k$th component $z_{j_k}$.
The convergence of the iterations (\ref{CD}) depends on the coordinate picking rule, namely how $j_k$ is chosen. Existing research works have analyzed the influence of different coordinate selection rules such as the cyclic rule and the random selection rule, on the convergence rate of the coordinate descent method. We choose the simplest variant using cyclic coordinate search to favor implementation simplicity. In fact, the cyclic implementation preserves the order of optimization variables with respect to the prediction horizon $t$, type (output $y$, input $u$, or input increment $\Delta u$), and component, so that reconstructing the coordinate that is currently optimized in is immediate. The implementation of one pass through all $n_z$ coordinates using cyclic CD is reported in Procedure~\ref{CCD}. In Procedure \ref{CCD}, the operator $\{s,\sigma\}=\operatorname{CCD}\{M,d,\sigma\}_{\underline{s}}^{\bar{s}}$ represents one pass iteration of the reverse cyclic CD method through all its $n_s$ coordinates $s_{1},\ldots,s_{n_s}$ for the box-constrained QP $\min _{s \in [\underline{s},\bar{s}]} \frac{1}{2} s^{\prime}Ms +s^{\prime} d$, that is to execute the following $n_s$ iterations
\begin{equation}\label{BoxCDM}
\begin{array}{l}
\text { for } i=1,\ldots,n_s \\
\quad \quad \hat{s}_i \leftarrow
\max(\underline{s}_i,\min(\bar{s}_i,s_i-\frac{1}{M_{i,i}}(M_{i,\cdot}s+d_i))\\
\quad \quad \sigma\leftarrow \sigma + (\hat{s}_i -s_i)^2\\
\quad \quad  s_i \leftarrow \hat{s}_i\\
\text { end }
\end{array}
\end{equation}

\floatname{algorithm}{Procedure}
\begin{algorithm}[t]
    \caption{Full pass of cyclic coordinate descent on all block variables}\label{CCD}
    \textbf{Input}: $\hat{\Lambda}=\{\hat{\lambda}_1,\ldots,\hat{\lambda}_{T}\}$, $\hat{\Gamma}=\{\hat{\gamma}_1,\ldots,\hat{\gamma}_{T}\}$, $Y =\{y_{1},  \cdots, y_{T}\}$, $U = \{u_{0},\cdots, u_{T-1}\}$, $\Delta U = \{\Delta u_{0},\cdots, \Delta u_{T-1}\}$; MPC settings $A(1),\ldots,A(n_a)$, $B(1),\ldots,B(n_b)$, $W^y$, $W^{\Delta u}$, $y_{\min}$, $y_{\max}$, $u_{\min}$, $u_{\max}$, $\Delta u_{\min}$, $\Delta u_{\max}$; parameter $\sigma, \rho>0$.
    \vspace*{.1cm}\hrule\vspace*{.1cm}
    \begin{enumerate}[label*=\arabic*., ref=\theenumi{}]
        \item $\sigma\leftarrow 0$;
        \item \textbf{for} $t=1,\ldots,T-1$ \textbf{do}
        \begin{enumerate}[label=\theenumi{}.\arabic*., ref=\theenumi{}.\arabic*]
            \item $j=\min(n_a, T-t)$;\label{CCD-jna}
            \item $\{y_{t}, \sigma\} \leftarrow {\operatorname{CCD}} \{\frac{1}{\rho}W^y+I+\sum_{i=1}^{j}A(i)^{\prime} A(i), e_t,\\ \sigma \}_{y_{\min}}^{y_{\max}}$;\label{CCD-yt}
            \item $j=\min(n_b, T-t+1)$;\label{CCD-jnb}
            \item $\{u_{t-1}, \sigma\} \leftarrow {\operatorname{CCD}} \{2I+\sum_{i=1}^{j}B(i)^{\prime} B(i), f_t,\\ \sigma \}_{u_{\min}}^{u_{\max}}$;\label{CCD-ut}
            \item $\{\Delta u_{t}, \sigma\} \leftarrow {\operatorname{CCD}} \{\frac{1}{\rho}W^{\Delta u}+I, g_t,\sigma \}_{\Delta u_{\min}}^{\Delta u_{\max}}$;\label{CCD-dut}
        \end{enumerate}
        \item $\{y_{T},\sigma\} \leftarrow {\operatorname{CCD}} \{\frac{1}{\rho}W^y+I, e_{T}, \sigma\}_{y_{\min}}^{y_{\max}}$;\label{CCD-yT}
        \item $\{u_{T-1},\sigma\} \leftarrow {\operatorname{CCD}} \{I+B(1)^{\prime}B(1), f_{T},\sigma \}_{u_{\min}}^{u_{\max}}$;\label{CCD-uT}
        \item $\{\Delta u_{T-1}, \sigma\} \leftarrow {\operatorname{CCD}} \{\frac{1}{\rho}W^{\Delta u}+I, g_{T}, \sigma \}_{\Delta u_{\min}}^{\Delta u_{\max}}$;\label{CCD-duT}
    \item \textbf{end}.
    \end{enumerate}
    \vspace*{.1cm}\hrule\vspace*{.1cm}
    \textbf{Output}: $Y$, $U$, $\Delta U$, $\Lambda$, $\Gamma$, $\sigma$.
\end{algorithm}
The quantities $e_t, f_t, g_t$ used in Procedure~\ref{CCD} are 
defined for $t=1,2,\cdots, T$ as follows:
\[\begin{aligned}
e_t&=-W^y r - (\lambda_t+\sum_{i=1}^{n_a}A(i)y_{t-i}+\sum_{i}^{n_b}B(i)u_{t-i})\\
&+ \sum_{n_i=1}^{\min(n_a,T-t)}A(n_i)^{\prime}(\lambda_{t+n_i}+\sum_{i\neq n_i}^{n_a}A(i)y_{t+n_i-i}\\
&+\sum_{i=1}^{n_b}B(i)u_{t+n_i-i}-y_{t+n_i})\\
f_t&=-(\gamma_{t}+u_{t-2}+\Delta u_{t-1})+(\gamma_{t+1}+\Delta u_{t}+u_{t})\\
&+ \sum_{n_i=1}^{\min(n_b,T-t+1)}B(n_i)^{\prime}(\lambda_{t+n_i}+\sum_{i=1}^{n_a}A(i)y_{t+n_i-i}\\
&+\sum_{i\neq n_i}^{n_b}B(i)u_{t+n_i-i}-y_{t+n_i})\\
g_t&=\gamma_{t}+u_{t-2}-u_{t-1}
\end{aligned}\]

\[\begin{aligned}
e_{T}&=-W^y r - (\lambda_{T}+\sum_{i=1}^{n_a}A(i)y_{T-i}+\sum_{i=1}^{n_b}B(i)u_{T-i})\\
f_{T}&=-(\gamma_{T}+u_{T-2}+\Delta u_{T-1})+B(1)^{\prime}(\lambda_{T}\\
&+\sum_{i=1}^{n_a}A(i)y_{T-i}+\sum_{i\neq1}^{n_b}B(i)u_{T-i}-y_{T})\\
g_{T}&=\gamma_{T}+u_{T-2}-u_{T-1}
\end{aligned}\]

which shows that they involve several matrix-vectors multiplications. It would greatly affect the computation efficiency since their computational cost is proportional to the product of the inner iterations and the outer iterations. 
To eliminate their explicit calculation, we propose here below an efficient coupling scheme between CD and AL that reduces the cost per iteration, without changing the rate of convergence of the algorithm.

\subsection{Efficient coupling scheme between CD and AL}
Our proposed efficient coupling scheme exploits the fact that CD only updates one coordinate each time, and the execution~\eqref{BoxCDM} of the operator $\operatorname{CCD}(\cdot)$ involves the next update of dual Lagrangian vectors. Here we take Step~\ref{CCD-yt} of Procedure~\ref{CCD} as an example, which has been modified from equation~\eqref{BoxCDM} to Procedure~\ref{CouplingCCD}. Note that the dual Lagrangian vectors used in Procedure~\ref{CouplingCCD} have been updated before Procedure~\ref{CouplingCCD}. The symbols $\{D^y_0,D^y_1,\ldots,D^y_{T-1}\}$ denote the diagonal elements of their Hessian matrices used in Step~\ref{CCD-yt} and~\ref{CCD-yT}
\begin{equation}
\begin{array}{l}
\text { for } t=1,\ldots,T-1 \\
\quad \quad j=\min(n_a, T-t);\\
\quad \quad D^y_t \leftarrow \diag\left(\frac{1}{\rho}W^y+I+\sum_{i=1}^{j}A(i)^{\prime}A(i)\right) \\
\text { end }\\
D^y_{T-1} \leftarrow \frac{1}{\rho}W^y+I
\end{array}
\end{equation}
To avoid repeating division operations, the values $\{\frac{1}{D^y_0},\frac{1}{D^y_1},\ldots,\frac{1}{D^y_{T-1}}\}$ are cached before the iterations start. The other steps involving the operator $\operatorname{CCD}(\cdot)$ in Procedure~\ref{CCD} follow the same idea.

\floatname{algorithm}{Procedure}
\begin{algorithm}[t]
    \caption{One pass of cyclic coordinate descent for Step \ref{CCD-yt} of Procedure \ref{CCD} after using efficient coupling scheme}\label{CouplingCCD}
    \textbf{Input}: $j=\min(n_a, T-t)$; $y_{t+1}$, $\lambda_t$, $\lambda_{t+1},\ldots,\lambda_{t+j}$; parameter $\rho>0$; update amount $\sigma\geq 0$.
    \vspace*{.1cm}\hrule\vspace*{.1cm}
    \begin{enumerate}[label*=\arabic*., ref=\theenumi{}]
        \item \textbf{for} $i=1,\ldots, n_y$ \textbf{do}
        \begin{enumerate}[label=\theenumi{}.\arabic*., ref=\theenumi{}.\arabic*]
            \item $s \leftarrow -\lambda_{t,i} + \sum_{n_i=1}^{j}A(n_i)_{:,i}^{\prime} \lambda_{t+n_i}$;
            \item $\theta \leftarrow  \left[y_{t,i}-\frac{\frac{1}{\rho}W^y_i(y_{t,i}-r_i) + s}{D^y_{t,i}}\right]_{y_{\min,i}}^{y_{\max,i}}$;
            \item $\Delta \leftarrow \theta - y_{t,i}$;
            \item $\sigma \leftarrow \sigma + \Delta^2$;
            \item $y_{t,i} \leftarrow \theta$;
            \item $\lambda_{t,i} \leftarrow \lambda_{t,i} + \Delta$;
            \item \textbf{for} $n_i=1,\ldots, j$ \textbf{do}
                \begin{enumerate}[label*=\arabic*., ref=\theenumi{}]
                \item $\lambda_{t+n_i} \leftarrow \lambda_{t+n_i} + \Delta \cdot A(n_i)_{:,i}$
                \end{enumerate}
        \end{enumerate}
    \item \textbf{end}.
    \end{enumerate}
    \vspace*{.1cm}\hrule\vspace*{.1cm}
    \textbf{Output}: $y_{t}$, $\lambda_t$, $\lambda_{t+1},\ldots,\lambda_{t+j}$, $\sigma$.
\end{algorithm}
\subsection{Algorithm}
Summarizing all the ingredients described in the previous sections, we obtain the construction-free ARX-based MPC Algorithm~\ref{A_CDAL_ARX}, which we call CDAL-ARX. 
Here, construction-free means that CDAL-ARX directly uses the ARX model coefficients without the need of constructing a QP problem explicitly. Note that the main update of the Lagrangian variables in Algorithm~\ref{A_CDAL_ARX} is placed early in Step~\ref{dual_update}, which is different from the original version of Algorithm 1
in~\cite{wu2021simple} because the CD method allows the use of our proposed efficient coupling scheme. The quantities $N_{\rm out}$ and $N_{\rm in}$ denote the maximum number of AL (outer-loop) and CD (inner-loop) iterations, respectively. The tolerances $\epsilon_{\rm out}$ and $\epsilon_{\rm in}$ define the stopping criteria of the outer and inner iterations, respectively. 

\floatname{algorithm}{Algorithm}
\begin{algorithm}
    \caption{Accelerated cyclic CDAL algorithm for ARX-based MPC}\label{A_CDAL_ARX}
    \textbf{Input}: primal/dual warm-start $Y = \{y_{1}, y_{2}, \cdots, y_{T}\}$, $U = \{u_{0}, u_{1}, \cdots, u_{T-1}\}$, $\Delta U =\{\Delta u_{0}, \Delta u_{1}, \cdots, \Delta u_{T-1}\}$, $\Lambda^{-1}$ $=$ $\Lambda^0$ $=$ $\{\lambda_{1}$, $\lambda_{2}$, $\cdots$, $\lambda_{T}\}$, $\Gamma^{-1}$ $=$ $\Gamma^0$ $=$ $\{\gamma_{1}$, $\gamma_{2}$, $\cdots$, $\gamma_{T}\}$; History input and output data $\{y_{0},y_{-1},\ldots,y_{1-n_a}\}$, $\{u_{-1},u_{-2},\ldots,u_{1-n_b}\}$; MPC settings $\{A(1), A(2),\ldots,A(n_a)$, $B(1),B(2),\ldots,B(n_b)$, $W^y$, $W^{\Delta u}$, $y_{\min}$, $y_{\max}$, $u_{\min}$, $u_{\max}$, $\Delta u_{\min}$, $\Delta u_{\max}\}$; Algorithm settings $\{\rho, N_{\rm out},N_{\rm in}\,\epsilon_{\rm out}, \epsilon_{\rm in}\}$
    \vspace*{.1cm}\hrule\vspace*{.1cm}
    \begin{enumerate}[label*=\arabic*., ref=\theenumi{}]
        \item $\alpha_1\leftarrow1$; $\hat{\Lambda}^{0}\leftarrow\Lambda^{0}$; $\hat{\Gamma}^{0}\leftarrow\Gamma^{0}$;
        \item \textbf{for} $k=1,2,\cdots, N_{\rm out}$ \textbf{do}
        \begin{enumerate}[label=\theenumi{}.\arabic*., ref=\theenumi{}.\arabic*]
            \item \textbf{for} $t=1,2,\cdots, T$ \textbf{do}\label{dual_update}
            \begin{enumerate}[label*=\arabic*., ref=\theenumi{}]
                \item $\lambda_{t}^{k}= \hat{\lambda}_{t}^{k-1}+
                (\sum_{i=1}^{n_{a}} A(i) y_{t-i}^{k}+
                \sum_{j=1}^{n_{b}} B(i) u_{t-i}^{k}-y_{t}^{k})$
                \item $\gamma_{t}^{k}=\hat{\gamma}_{t}^{k-1}+(u_{t-2}^{k}+\Delta u_{t-1}^{k}-u_{t-1}^{k})$;
            \end{enumerate}
            \item \textbf{for}  $k_{in}=1,2,\cdots, N_{\rm in}$  \textbf{do}         
            \begin{enumerate}[label=\theenumii{}.\arabic*., ref=\theenumii{}.\arabic*]
                \item $(Y, U, \Delta U, \sigma) \leftarrow$  Procedure \ref{CCD} with use of Procedure \ref{CouplingCCD};
                \item \textbf{if} $\sigma \leq \epsilon_{\rm in}$ \textbf{break} the loop;
            \end{enumerate}            
            \item \textbf{if} $\|\Lambda^{k}-\hat{\Lambda}^{k-1}\|_2^2\leq \epsilon_{\rm out}$ \textbf{stop};
            \item $\alpha_{k+1} \leftarrow \frac{1+\sqrt{1+4\alpha_k^2}}{2}$;
            \item $\hat{\Lambda}^{k} \leftarrow \Lambda^k + \frac{\alpha_k-1}{\alpha_{k+1}} (\Lambda^k-\Lambda^{k-1})$;
        \end{enumerate}
    \item \textbf{end}.
    \end{enumerate}
    \vspace*{.1cm}\hrule\vspace*{.1cm}
    ~~\textbf{Output}: $Y,U,\Delta U,\Lambda,\Gamma$
\end{algorithm}

\section{Numerical examples}\label{sec:example}
In this section, we test our proposed ARX-based MPC algorithm against other MPC solvers, which rely on condensed or sparse MPC-to-QP construction, respectively. The best choice between condensed and sparse QP forms mainly depends on the number of outputs $n_y$, control inputs $n_u$, and the length of the prediction horizon $T$~\cite{kouzoupis2015towards}. For numerical comparisons with our ARX-based MPC algorithm, this paper considers both condensed and sparse MPC-to-QP constructions, which are then solved by the qpOASES~\cite{ferreau2014qpoases} and OSQP~\cite{stellato2020osqp}, respectively. The reported comparison simulation results were obtained on a MacBook Pro with a 2.7~GHz 4-core Intel Core i7 and 16GB RAM. Algorithm~\ref{A_CDAL_ARX}, qpOASES v3.2 and OSQP v0.6.2 are all executed in MATLAB R2020a
via their C-mex implementations.

\subsection{Problem descriptions}
\begin{enumerate}
    \item Time-varying ARX model example: one notable feature of ARX models is their ease to be updated at runtime, which makes them particularly appealing when the system dynamics cannot be well captured by a single linear time-invariant model. Our CDAL-ARX algorithm can take advantage of its construction-free feature to avoid the computation cost of the online construction step. We tested CDAL-ARX on randomly-generated two-input-two-output ARX models with order $n_a=4$ and $n_b=4$ and time-varying system matrices. For demonstration purposes, here below we report one instance of them, whose ARX coefficient matrices $A^t(1),\ldots,A^t(4)$, $B^t(1),\ldots,B^t(4)$ at time $t$ are given by
    \begin{equation}\label{eqn_time_varying_ARX}
    \begin{array}{l}
    A(i)^{t}=A(i)+0.1M^t, i=1,\ldots 4 \\
    B(i)^{t}=B(i)+0.1M^t, i=1,\ldots 4
    \end{array}
    \end{equation}
    where $A(1)\!=\!\scriptsize\left[\begin{array}{cc}
    0.9&0.1\\
    0.1&0.9
    \end{array}\right]\!,
    A(2)\!=\!\scriptsize\left[\begin{array}{cc}
    0.7&0.1\\
    0.1&0.7
    \end{array}\right],A(3)$\\
    $\!=\!\scriptsize\left[\begin{array}{cc}
        0.5&0.1\\
        0.1&0.5
        \end{array}\right],      
    A(4)\!=\!{\scriptsize\left[\begin{array}{cc}
    0.3&0.1\\
    0.1&0.3
    \end{array}\right]},B(1)\!=\!$\\
    $\scriptsize\left[\begin{array}{cc}
    1&0.5\\
    0.5&1
    \end{array}\right],B(2)\!=\!\scriptsize\left[\begin{array}{cc}
    0.8&0.4\\
    0.4&0.8
    \end{array}\right],B(3)\!=\!\scriptsize\left[\begin{array}{cc}
    0.6&0.3\\
    0.3&0.6
    \end{array}\right]$\\
    $,B(4)\!=\!\scriptsize\left[\begin{array}{cc}
    0.4&0.2\\
    0.2&0.4
    \end{array}\right],M^t\!=\!\scriptsize\left[\begin{array}{cc}
    \sin(\frac{t}{10})& \cos(\frac{t}{10})\\[1em]
    \cos(\frac{t}{10})&\sin(\frac{t}{10})
    \end{array}\right].$\\


    \item \textit{DNN-based LPV-ARX model example}: We tested on CDAL-ARX on randomly-generated two-input-two-output quasi-LPV-ARX models of larger order $n_a=6$ and $n_b=6$, whose coefficient matrices are piecewise affine (PWA) maps of the scheduling vector $w_{t-1}$
    \begin{equation}\label{DNN_LPV_ARX}
    \left[\begin{array}{c}
        y_t(1)  \\
        y_t(2) 
    \end{array}\right]
    = \left[\begin{array}{c}
         \mathcal{N}_1(w_{t-1})^{\prime}  \\
         \mathcal{N}_2(w_{t-1})^{\prime}
    \end{array}
    \right] x_{t-1}
    \end{equation}
    where $x_{t-1}\!=\! \left[y_{t-1}^{\prime} ,\ldots,y_{t-6}^{\prime},u_{t-1}^{\prime},\ldots,u_{t-6}^{\prime}\right]^{\prime} \in \mathbb{R}^{24}$, $w_{t-1}=\left[y_{t-1}^{\prime},\ldots,y_{t-6}^{\prime},u_{t-2}^{\prime},\ldots,u_{t-6}^{\prime}\right]^{\prime} \in \mathbb{R}^{22}$, and $\mathcal{N}_1$, $\mathcal{N}_2\in \mathbb{R}^{22} \to\mathbb{R}^{24}$ are deep feedforward neural networks with three layers and ReLU activation function, namely $\mathcal{N}_1(w_{t-1})\!=\!W_{1,3}\max(0,W_{1,2}\max(0,W_{1,1}v_{t-1}+b_{1,1})$
    $+b_{1,2})+b_{1,3}$, $\mathcal{N}_2(w_{t-1})\!=\!W_{2,3}\max(0,W_{2,2}\max(0,$\\
    $W_{2,1}v_{t-1}+b_{2,1})+b_{2,2})+b_{2,3}$. Here we choose the number of neurons in each hidden layer as three times the number of inputs according to \cite{serra2018bounding}, that is, $W_{1,1}$ and $ W_{2,1}\in\mathbb{R}^{66\times22}$, $b_{1,1}$ and $b_{2,1}\in\mathbb{R}^{66}$, $W_{1,2}$ and $W_{2,2}\in\mathbb{R}^{66\times66}$, $b_{1,2}$ and $b_{2,2}\in\mathbb{R}^{66}$, $W_{1,3}$ and $W_{2,3}\in\mathbb{R}^{24\times66}$, $b_{1,3}$ and $b_{2,3}\in\mathbb{R}^{24}$. For demonstration purposes, we define $b_{1,3},b_{2,3}$ by collecting the coefficients defining $A(1),\ldots,A(4),A(4),A(4),B(1),\ldots,B(4),B(4),B(4)$ as in~\eqref{eqn_time_varying_ARX}; the remaining network parameters are randomly generated uniformly between 0 and 0.1. At each time $t$, the linear model consumed by our CDAL-ARX algorithm is given by evaluating the deep ReLU networks as in~\eqref{DNN_LPV_ARX}.
\end{enumerate}

\begin{figure*}[htbp]
\centering 
\subfigure[Time-varying ARX model]{
\begin{minipage}{8cm}\label{fig1}
\centering
\includegraphics[scale=0.43]{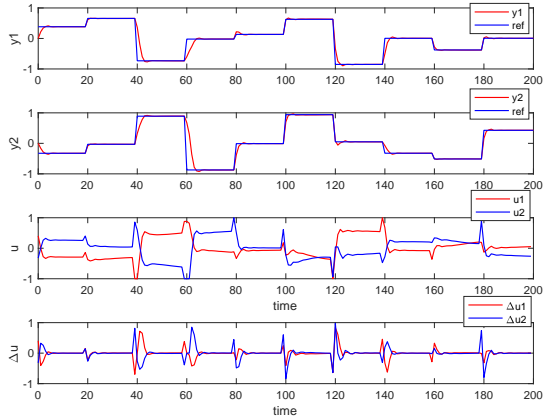}
\end{minipage}
}\subfigure[DNN-based LPV-ARX model]{
\begin{minipage}{8cm}\label{fig2}
\centering 
\includegraphics[scale=0.43]{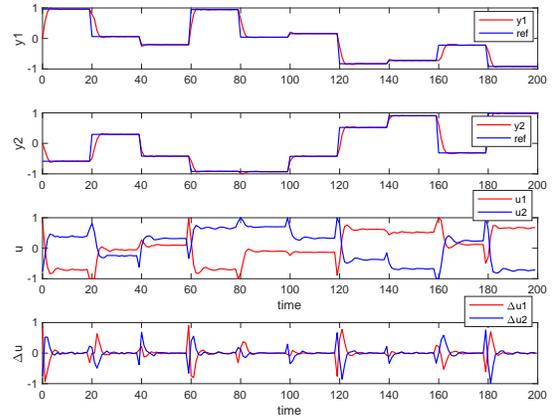}
\end{minipage}
}
\caption{Closed-loop tracking results}
\end{figure*}

In both examples, we use the same MPC parameters $W_y=I$, $W_{\Delta u}=0.1I$, $[y_{\rm min},y_{\rm max}]=[-1,1]$, $[u_{\rm min},u_{\rm max}]=[-1,1]$, $[\Delta u_{\rm min},\Delta u_{\rm max}]=[-1,1]$. Different prediction horizon lengths $T$ are used to investigate numerical performance, namely $T=10$, $20$, and $30$. The initial conditions are $y_{-3}=y_{-2}=y_{-1}=y_{0}=[0\ 0]^{\prime}$, and $u_{-3}=u_{-2}=u_{-1}=[0\ 0]^{\prime}$. In the two examples, the closed-loop simulation is run over $200$ sampling steps, and the desired references for $y_1$ and $y_2$ are randomly changed every 20 steps. Warm-start used in  all solvers (qpOASES, OSQP, CDAL-ARX). 
We keep default solver settings in both qpOASES and OSQP, so that they produce solutions of similar precision, that is measured in terms of Euclidean distance (since qpOASES belongs to the class of active-set methods, in principle it always provides a high-precision solution at termination, so its solution quality cannot be tuned as easily as in the case of ADMM). For a fair comparison, in the two examples we set
$\epsilon_{in}=10^{-6}$ and $\epsilon_{out}=10^{-6}$ under $\rho=1$ to define
the stopping criteria of our CDAL-ARX solver, so to obtain closed-loop control sequences with similar precision. In both examples, the generated closed-loop simulation results are almost indistinguishable, see Figures~\ref{fig1} and~\ref{fig2}, respectively, which show good tracking performance and no violation in input and output constraints.

Using the qpOASES and OSQP solvers require the online construction of the QP problem, whose computation time must be counted in the total time. Table~\ref{Tab1} lists the solution time of CDAL-ARX and lists the construction and solution time when using qpOASES (condensed construction) and OSQP (sparse construction). From Table~\ref{Tab1} it can be noticed that CDAL-ARX is always solving the MPC problem in a smaller CPU time, when compared to the sum of construction and solution time of qpOASES and OSQP. Moreover, as the prediction horizon increases, qpOASES and OSQP may fail to solve the problem due to the ill-conditioning issue. Note also that the computation time of CDAL-ARX is often shorter than the pure solution time of qpOASES and OSQP (i.e., not counting the construction time), which seems to indicate that the reported speed-ups are due to both adopting the proposed augmented Lagrangian method and avoiding the use of state-space models.

\begin{table}[!htbp]
\caption{Computation time (ms) of CDAL-ARX and comparison with other solvers}
\centering
\begin{threeparttable}
\begin{tabular}{ccccc}
\toprule
Examples & T &  CDAL-ARX  &  qpOASES & OSQP\\
&  &  avg, max & avg, max & avg, max\\
\midrule
\multirow{6}*{Time-varying ARX} & 10 & 0.14, 1.5 & 0.42\tnote{*}, 2.8\tnote{*} & 0.41\tnote{*}, 2.9\tnote{*}\\
~ & & & 0.08\tnote{\dag}, 1.4\tnote{\dag} & 0.18\tnote{\dag}, 0.92\tnote{\dag}  \\
~ & 20 & 0.25, 2.8 & 1.2\tnote{*}, 6.3\tnote{*} & 1.0\tnote{*}, 3.8\tnote{*}\\
~ & & & 0.18\tnote{\dag}, 4.4\tnote{\dag} & 1.9\tnote{\dag}, 17\tnote{\dag}  \\
~ & 30 & 0.36, 3.6 & 2.6\tnote{*}, 10.2\tnote{*} & 2.4\tnote{*}, 8.1\tnote{*}\\
~ & & & fail & 22\tnote{\dag}, 48\tnote{\dag}  \\
\midrule
\multirow{6}*{DNN-based LPV-ARX} & 10 & 0.51, 2.5 & 0.46\tnote{*}, 3.2\tnote{*} & 0.42\tnote{*}, 3.6\tnote{*}\\
~ & & & 0.57\tnote{\dag}, 3.9\tnote{\dag} & 1.1\tnote{\dag}, 14\tnote{\dag}  \\
~ & 20 & 1.2, 4.6 & 1.2\tnote{*}, 5.5\tnote{*} & 0.97\tnote{*}, 4.5\tnote{*}\\
~ & & & fail & 16\tnote{\dag}, 32\tnote{\dag}  \\
~ & 30 & 2.0, 7.3 & 3.1\tnote{*}, 10.8\tnote{*} & 2.8\tnote{*}, 8.9\tnote{*}\\
~ & & & fail & fail  \\
\bottomrule
\end{tabular}
\begin{tablenotes}
    \footnotesize
    \item \tnote{*}construction time, \tnote{\dag}solution time. For qpOASES and OSQP the time to evaluate the MPC law is the sum of construction and solution time.
\end{tablenotes}
\end{threeparttable}
\label{Tab1}
\end{table}

\section{Conclusion}\label{sec:conclusion}
This paper has introduced a solution algorithm for solving MPC problems based on ARX models that avoids constructing the associated QP problem explicitly. Due to its matrix-free and library-free features, the proposed CDAL-ARX algorithm can be useful in adaptive embedded linear MPC applications based on ARX models, especially when combined with a fast and robust recursive linear identification method. Future research will address extending the method to handle soft output constraints, so to relax Assumption~\ref{assump1}.

\bibliographystyle{unsrt}
\bibliography{ref} 

\end{document}